\documentclass[english,11pt]{article}
\usepackage{amsfonts}
\usepackage[latin1]{inputenc}
\usepackage{amsmath,amsthm,amssymb}
\usepackage{graphicx}
\usepackage{color}

\def\neweq#1{\begin{equation}\label{#1}}
\def\endeq{\end{equation}}

\newtheorem{theorem}{Theorem}[section]

\newtheorem{remark}{Remark}[section]
\input{tcilatex}
\begin{document}

\title{\textbf{An existence result for a quasilinear system with gradient
term under the Keller-Osserman conditions}}
\author{Dragos-Patru Covei\break \\
{\small \ Constantin Brancusi University of \ Tg-Jiu, Calea Eroilor Nr.30,
Targu-Jiu, Gorj, Romania }\\
{\small E-mail: patrucovei\texttt{@yahoo.com}}}
\date{}
\maketitle

\begin{abstract}
We deal with existence of entire solutions for the quasilinear elliptic
system of this type $\Delta _{p}u_{i}+h_{i}\left( \left\vert x\right\vert
\right) \left\vert \nabla u_{i}\right\vert ^{p-1}=a_{i}\left( \left\vert
x\right\vert \right) f_{i}\left( u_{1},u_{2}\right) $ on $\mathbb{R}^{N}$ ($%
N\geq 3$, $i=1,2$) where $N-1\geq p>1$, $\Delta _{p}$ is the p-Laplacian
operator and $h_{i}$, $a_{i}$, $f_{i}$ are suitable functions. The results
of this paper supplement the existing results in the literature and improve
those obtained by Xinguang Zhang and Lishan Liu, The existence and
nonexistence of entire positive solutions of semilinear elliptic systems
with gradient term, Journal of Mathematical Analysis and Applications,
Volume 371, Issue 1, 1 November 2010, Pages 300-308).
\end{abstract}

\baselineskip16pt \renewcommand{\theequation}{\arabic{section}.%
\arabic{equation}} \catcode`@=11 \@addtoreset{equation}{section} \catcode%
`@=12

\textbf{2000 Mathematics Subject Classification}:
35J60;35J62;35J66;35J92;58J10;58J20.

\textbf{Key words}: Entire solution; Large solution; Elliptic system.

\section{Introduction}

In the present paper we establish a new result concerning the existence of
solutions for the quasilinear elliptic system 
\begin{equation}
\left\{ 
\begin{array}{l}
\Delta _{p}u_{1}\left( r\right) +h_{1}\left( r\right) \left\vert \nabla
u_{1}\left( r\right) \right\vert ^{p-1}=a_{1}\left( r\right) f_{1}\left(
u_{1}\left( r\right) ,u_{2}\left( r\right) \right) \text{ }, \\ 
\Delta _{p}u_{2}\left( r\right) +h_{2}\left( r\right) \left\vert \nabla
u_{2}\left( r\right) \right\vert ^{p-1}=a_{2}\left( r\right) f_{2}\left(
u_{1}\left( r\right) ,u_{2}\left( r\right) \right) \text{ },%
\end{array}%
\right. \text{ }  \label{11}
\end{equation}%
where $r:=\left\vert x\right\vert $ for $x\in \mathbb{R}^{N}$ ($N-1\geq p>1$%
) is the Euclidean norm, $\Delta _{p}$ is the so called p-Laplace operator
defined by $\Delta _{p}u:=\func{div}\left( \left\vert \nabla u\right\vert
^{p-2}\nabla u\right) $. It will be assumed throughout this paper that $%
a_{j} $, $h_{j}$ ($j=1,2$) are nonnegative nontrivial $C\left( \mathbb{R}%
^{N}\right) $ functions, while $f_{j}:\left[ 0,\infty \right)
^{2}\rightarrow \left[ 0,\infty \right) $ ($j=1,2$) are continuous and
nondecreasing functions in each variable and verify $f_{j}\left(
s_{1},s_{2}\right) >0$ whenever $s_{i}>0$ for some $i=1,2$ together with the
"Keller-Osserman type" condition

\begin{equation}
\quad I\left( \infty \right) :=\lim_{r\rightarrow \infty }I\left( r\right)
=\infty  \label{KO}
\end{equation}%
where $I\left( r\right) :=\int_{a}^{r}[F\left( s\right) ]^{-1/2}ds$ for $%
r\geq a>0$, $F\left( s\right) :=\int_{0}^{s}\underset{i=1}{\overset{2}{%
\Sigma }}f_{i}\left( t,t\right) dt$.

For a single equation of the form $\Delta u=f\left( u\right) $ where $%
f\left( u\right) $\ is positive, real continuous function defined for all
real $u$ and nondecreasing the existence of entire large solutions is
equivalent to a condition on $f$ known as the Keller--Osserman condition 
\begin{equation}
\int_{u_{0}}^{\infty }\left( \int_{0}^{t}f\left( s\right) ds\right)
^{-1/2}dt=\infty \text{ for \ }u_{0}>0,  \label{KO1}
\end{equation}%
(see \cite{K}, \cite{O}). In particular, Keller and Osserman prove that a
necessary and sufficient condition for the considered problem to have an
entire solution is that $f$ satisfies (\ref{KO1}). Such a solution will
necessarily satisfies $\lim_{\left\vert x\right\vert \rightarrow \infty
}u\left( x\right) =\infty $ and hence be a large solution. Moreover, Keller
applied the results to electrohydrodynamics, namely to the problem of the
equilibrium of a charged gas in a conducting container, see \cite{K1}.

There is by now a broad literature regarding the study of solutions for (\ref%
{11}). Basic results in the study of solutions for such problems have been
obtained in the last few decades in the works of \cite%
{BEC4,CD2,KN,JM1C4,JM1,NI,LZZ} and their references. We comment below on a
few further results closer to our interests in the present article.

Regarding (\ref{11}), Zhang and Liu \cite{LZZ} studied the existence of
entire large positive solutions of the system%
\begin{equation*}
\text{ }\left\{ 
\begin{array}{l}
\Delta u_{1}+\left\vert \nabla u_{1}\right\vert =a_{1}\left( r\right)
f_{1}\left( u_{1},u_{2}\right) \text{,} \\ 
\Delta u_{2}+\left\vert \nabla u_{2}\right\vert =a_{2}\left( r\right)
f_{2}\left( u_{1},u_{2}\right) \text{. }%
\end{array}%
\right.
\end{equation*}%
In \cite{LZZ}, the authors imposed on $a_{1}$, $a_{2}$, $f_{1}$ and $f_{2}$
satisfying the above conditions and instead of the Keller-Osserman condition
the following 
\begin{equation}
\int_{a}^{\infty }\frac{ds}{f_{1}\left( s,s\right) +f_{2}\left( s,s\right) }%
=\infty \text{ for }r\geq a>0.  \label{LZZ}
\end{equation}%
Obviously, (\ref{LZZ}) implies (\ref{KO}).

Finally, we note that the study of large solutions for (\ref{11}) when the
integral in (\ref{KO}) is finite has been the subject of the article \cite{W}%
.

Motivated by papers \cite{W} and \cite{LZZ} we are interested in another
type of nonlinearity $f_{i}$ ($i=1,2$) in order to obtain the existence of
entire large/bounded positive solutions of (\ref{11}).

The main reult of this article is:

\begin{theorem}
\label{1}Under the above hypotheses \textit{there are infinitely many
positive entire radial solutions of system (\ref{11}). Suppose furthermore
that }$r^{\frac{p\left( N-1\right) }{p-1}}\underset{j=1}{\overset{2}{\Sigma }%
}e^{\frac{p}{p-1}\int_{0}^{r}h_{j}\left( t\right) dt}a_{j}\left( r\right) $%
\textit{\ is nondecreasing for large }$r$. \textit{Then the solutions}:

i)\quad are bounded if \textit{there exists a positive number }$\varepsilon $%
\textit{\ such that }%
\begin{equation}
\int_{0}^{\infty }t^{1+\varepsilon }\left( \underset{j=1}{\overset{2}{\Sigma 
}}e^{\frac{p}{p-1}\int_{0}^{t}h_{j}(t)dt}a_{j}\left( t\right) \right)
^{2/p}dt<\infty ,\text{ }  \label{5}
\end{equation}%
ii)\quad are large if \textit{\ }%
\begin{equation}
\int_{0}^{\infty }\left( \frac{e^{-\int_{0}^{t}h_{j}\left( s\right) ds}}{%
t^{N-1}}\int_{0}^{t}s^{N-1}e^{\int_{0}^{s}h_{j}\left( t\right)
dt}a_{j}\left( s\right) ds\right) ^{1/\left( p-1\right) }dt=\infty \text{
for all }j=1,2  \label{12}
\end{equation}%
hold.
\end{theorem}

Our main result are new, because no solutions have been detected yet for the
system of the form (\ref{11}) under the Keller-Osserman conditions (\ref{KO}%
). We mention that we can prove similar results for $f_{1}$ and $f_{2}$
being non-monotonic as in \cite{K}, \cite{W}. Since in this case the proof
is as for the monotone case we omit it.

\section{Proof of the Theorem \protect\ref{1} \label{ss1}}

We start by showing that (\ref{11}) has positive radial solutions. The proof
is inspired by \cite{CD2} with some new ideas. Note that radial solutions of
(\ref{11}) are radial solutions of the system 
\begin{equation}
\left\{ 
\begin{array}{c}
\left( p-1\right) u_{1}^{\prime }\left( r\right) ^{p-2}u_{1}^{\prime \prime
}+\frac{N-1}{r}u_{1}^{\prime }\left( r\right) ^{p-1}+h_{1}\left( r\right)
\left\vert u_{1}^{\prime }\left( r\right) \right\vert ^{p-1}=a_{1}\left(
r\right) f_{1}\left( u_{1}\left( r\right) ,u_{2}\left( r\right) \right) , \\ 
\left( p-1\right) u_{2}^{\prime }\left( r\right) ^{p-2}u_{2}^{\prime \prime
}+\frac{N-1}{r}u_{2}^{\prime }\left( r\right) ^{p-1}+h_{2}\left( r\right)
\left\vert u_{2}^{\prime }\left( r\right) \right\vert ^{p-1}=a_{2}\left(
r\right) f_{2}\left( u_{1}\left( r\right) ,u_{2}\left( r\right) \right) ,%
\end{array}%
\right.  \label{66}
\end{equation}%
where we can assume in the next that $u_{i}^{\prime }\left( r\right) \geq 0$
($i=1,2$).

First we see that radial solutions of (\ref{66}) are any positive solutions $%
\left( u_{1},u_{2}\right) $ of the integral equations%
\begin{equation}
\left\{ 
\begin{array}{c}
u_{1}\left( r\right) =\frac{b}{2}+\int_{0}^{r}\left( \frac{%
e^{-\int_{0}^{t}h_{1}\left( s\right) ds}}{t^{N-1}}\int_{0}^{t}s^{N-1}e^{%
\int_{0}^{s}h_{1}\left( s\right) dt}a_{1}\left( s\right) f_{1}\left(
u_{1}\left( s\right) ,u_{2}\left( s\right) \right) ds\right) ^{1/\left(
p-1\right) }dt, \\ 
u_{2}\left( r\right) =\frac{b}{2}+\int_{0}^{r}\left( \frac{%
e^{-\int_{0}^{t}h_{2}\left( s\right) ds}}{t^{N-1}}\int_{0}^{t}s^{N-1}e^{%
\int_{0}^{s}h_{2}\left( s\right) dt}a_{2}\left( s\right) f_{2}\left(
u_{1}\left( s\right) ,u_{2}\left( s\right) \right) ds\right) ^{1/\left(
p-1\right) }dt,%
\end{array}%
\right.  \label{non}
\end{equation}%
where $b\geq a>0$. Our idea is to regard this as an operator equation 
\begin{equation*}
S\left( u_{1}(r),u_{2}\left( r\right) \right) =\left( u_{1}(r),u_{2}\left(
r\right) \right)
\end{equation*}%
with 
\begin{equation*}
S:C\left[ 0,\infty \right) \times C\left[ 0,\infty \right) \rightarrow C%
\left[ 0,\infty \right) \times C\left[ 0,\infty \right) \text{ }
\end{equation*}%
defined by%
\begin{equation}
S\left( u_{1}(r),u_{2}\left( r\right) \right) =\left( 
\begin{array}{c}
\frac{b}{2}+\int_{0}^{r}\left( \frac{e^{-\int_{0}^{t}h_{1}\left( s\right) ds}%
}{t^{N-1}}\int_{0}^{t}s^{N-1}e^{\int_{0}^{s}h_{1}\left( s\right)
dt}a_{1}\left( s\right) f_{1}\left( u_{1}\left( s\right) ,u_{2}\left(
s\right) \right) ds\right) ^{1/\left( p-1\right) }dt, \\ 
\frac{b}{2}+\int_{0}^{r}\left( \frac{e^{-\int_{0}^{t}h_{2}\left( s\right) ds}%
}{t^{N-1}}\int_{0}^{t}s^{N-1}e^{\int_{0}^{s}h_{2}\left( s\right)
dt}a_{2}\left( s\right) f_{2}\left( u_{1}\left( s\right) ,u_{2}\left(
s\right) \right) ds\right) ^{1/\left( p-1\right) }dt%
\end{array}%
\right) ^{T}  \label{op}
\end{equation}%
where $u_{1}(0)=\frac{b}{2}$ and $u_{2}(0)=\frac{b}{2}$ are the central
values for the system. The integration in this operator implies that a fixed
point $\left( u_{1},u_{2}\right) \in C\left[ 0,\infty \right) \times C\left[
0,\infty \right) $ is in fact in the space $C^{1}\left[ 0,\infty \right)
\times C^{1}\left[ 0,\infty \right) $. Then a solution of (\ref{66}) will be
obtained as a fixed point of the operator (\ref{op}). To establish a
solution to this operator, we use successive approximation. We define,
recursively, sequences $\left\{ u_{i}^{k}\right\} _{i=\overline{1,2}}^{k\geq
1}$ on $\left[ 0,\infty \right) $ by 
\begin{equation*}
u_{1}^{0}=u_{2}^{0}=\frac{b}{2}\text{ for all }r\geq 0
\end{equation*}%
and 
\begin{eqnarray*}
\left( u_{1}^{k}\left( r\right) ,u_{2}^{k}\left( r\right) \right) &=&S\left(
u_{1}^{k-1}\left( r\right) ,u_{2}^{k-1}\left( r\right) \right) \\
&=&\left( 
\begin{array}{c}
\frac{b}{2}+\int_{0}^{r}\left( \frac{e^{-\int_{0}^{t}h_{1}\left( s\right) ds}%
}{t^{N-1}}\int_{0}^{t}s^{N-1}e^{\int_{0}^{s}h_{1}\left( s\right)
dt}a_{1}\left( s\right) f_{1}\left( u_{1}^{k-1}\left( s\right)
,u_{2}^{k-1}\left( s\right) \right) ds\right) ^{1/\left( p-1\right) }dt \\ 
\frac{b}{2}+\int_{0}^{r}\left( \frac{e^{-\int_{0}^{t}h_{2}\left( s\right) ds}%
}{t^{N-1}}\int_{0}^{t}s^{N-1}e^{\int_{0}^{s}h_{2}\left( s\right)
dt}a_{2}\left( s\right) f_{2}\left( u_{1}^{k-1}\left( s\right)
,u_{2}^{k-1}\left( s\right) \right) ds\right) ^{1/\left( p-1\right) }dt%
\end{array}%
\right) ^{T}.
\end{eqnarray*}%
We remark that, for all $r\geq 0,$ $i=1,2$ and $k\in N$ 
\begin{equation*}
u_{i}^{k}\left( r\right) \geq \frac{b}{2}\text{,}
\end{equation*}%
and that $\left\{ u_{i}^{k}\right\} _{i=1,2}^{k\geq 1}$ is an increasing
sequence of nonnegative, non-decreasing functions.

We note that $\left\{ u_{i}^{k}\right\} _{i=1,2}^{k\geq 1}$ satisfy%
\begin{equation}
\left\{ 
\begin{array}{c}
\left( p-1\right) \left[ \left( u_{1}^{k}\right) ^{\prime }\right]
^{p-2}\left( u_{1}^{k}\right) ^{\prime \prime }+\left( \frac{N-1}{r}%
+h_{1}\left( r\right) \right) \left[ \left( u_{1}^{k}\right) ^{\prime }%
\right] ^{p-1}=a_{1}\left( r\right) f_{1}\left( u_{1}^{k-1}\left( r\right)
,u_{2}^{k-1}\left( r\right) \right) , \\ 
\left( p-1\right) \left[ \left( u_{2}^{k}\right) ^{\prime }\right]
^{p-2}\left( u_{2}^{k}\right) ^{\prime \prime }+\left( \frac{N-1}{r}%
+h_{1}\left( r\right) \right) \left[ \left( u_{2}^{k}\right) ^{\prime }%
\right] ^{p-1}=a_{2}\left( r\right) f_{2}\left( u_{1}^{k-1}\left( r\right)
,u_{2}^{k-1}\left( r\right) \right) .%
\end{array}%
\right.  \label{sis1}
\end{equation}%
Using the monotonicity of $\left\{ u_{i}^{k}\right\} _{i=1,2}^{k\geq 1}$
yields%
\begin{eqnarray}
a_{1}\left( r\right) f_{1}\left( u_{1}^{k-1}\left( r\right)
,u_{2}^{k-1}\left( r\right) \right) &\leq &a_{1}\left( r\right) f_{1}\left(
u_{1}^{k},u_{2}^{k}\right) \leq a_{1}\left( r\right) \overset{2}{\underset{%
i=1}{\Sigma }}f_{i}\left( \overset{2}{\underset{i=1}{\Sigma }}u_{i}^{k},%
\overset{2}{\underset{i=1}{\Sigma }}u_{i}^{k}\right) ,  \notag \\
&&  \label{8} \\
a_{2}\left( r\right) f_{2}\left( u_{1}^{k-1}\left( r\right)
,u_{2}^{k-1}\left( r\right) \right) &\leq &a_{2}\left( r\right) f_{2}\left(
u_{1}^{k},u_{2}^{k}\right) \leq a_{2}\left( r\right) \overset{2}{\underset{%
i=1}{\Sigma }}f_{i}\left( \overset{2}{\underset{i=1}{\Sigma }}u_{i}^{k},%
\overset{2}{\underset{i=1}{\Sigma }}u_{i}^{k}\right) ,  \notag
\end{eqnarray}%
and, so%
\begin{equation}
\left\{ 
\begin{array}{c}
\left( p-1\right) \left[ \left( u_{1}^{k}\left( r\right) \right) ^{\prime }%
\right] ^{p-1}\left( u_{1}^{k}\right) ^{\prime \prime }+\left( \frac{N-1}{r}%
+h_{1}\left( r\right) \right) \left[ \left( u_{1}^{k}\left( r\right) \right)
^{\prime }\right] ^{p}\leq a_{1}\left( r\right) \overset{2}{\underset{i=1}{%
\Sigma }}f_{i}\left( \overset{2}{\underset{i=1}{\Sigma }}u_{i}^{k},\overset{2%
}{\underset{i=1}{\Sigma }}u_{i}^{k}\right) \left( u_{1}^{k}\left( r\right)
\right) ^{\prime }, \\ 
\left( p-1\right) \left[ \left( u_{2}^{k}\left( r\right) \right) ^{\prime }%
\right] ^{p-1}\left( u_{2}^{k}\right) ^{\prime \prime }+\left( \frac{N-1}{r}%
+h_{2}\left( r\right) \right) \left[ \left( u_{2}^{k}\left( r\right) \right)
^{\prime }\right] ^{p}\leq a_{2}\left( r\right) \overset{2}{\underset{i=1}{%
\Sigma }}f_{i}\left( \overset{2}{\underset{i=1}{\Sigma }}u_{i}^{k},\overset{2%
}{\underset{i=1}{\Sigma }}u_{i}^{k}\right) \left( u_{2}^{k}\left( r\right)
\right) ^{\prime },%
\end{array}%
\right.  \label{88}
\end{equation}%
which implies that%
\begin{equation}
\left\{ 
\begin{array}{c}
\left( p-1\right) \left[ \left( u_{1}^{k}\left( r\right) \right) ^{\prime }%
\right] ^{p-1}\left( u_{1}^{k}\right) ^{\prime \prime }+\left( \frac{N-1}{r}%
+h_{1}\left( r\right) \right) \left[ \left( u_{1}^{k}\left( r\right) \right)
^{\prime }\right] ^{p}\leq a_{1}\left( r\right) \overset{2}{\underset{i=1}{%
\Sigma }}f_{i}\left( \overset{2}{\underset{i=1}{\Sigma }}u_{i}^{k},\overset{2%
}{\underset{i=1}{\Sigma }}u_{i}^{k}\right) \left( \overset{2}{\underset{i=1}{%
\Sigma }}u_{i}^{k}\left( r\right) \right) ^{\prime }, \\ 
\left( p-1\right) \left[ \left( u_{2}^{k}\left( r\right) \right) ^{\prime }%
\right] ^{p-1}\left( u_{2}^{k}\right) ^{\prime \prime }+\left( \frac{N-1}{r}%
+h_{2}\left( r\right) \right) \left[ \left( u_{2}^{k}\left( r\right) \right)
^{\prime }\right] ^{p}\leq a_{2}\left( r\right) \overset{2}{\underset{i=1}{%
\Sigma }}f_{i}\left( \overset{2}{\underset{i=1}{\Sigma }}u_{i}^{k},\overset{2%
}{\underset{i=1}{\Sigma }}u_{i}^{k}\right) \left( \overset{2}{\underset{i=1}{%
\Sigma }}u_{i}^{k}\left( r\right) \right) ^{\prime }.%
\end{array}%
\right.  \label{sum}
\end{equation}%
Let 
\begin{equation*}
a_{i}^{R}=\max \{a_{i}\left( r\right) :0\leq r\leq R\}\text{, }i=1,2\text{.}
\end{equation*}%
We prove that $u_{i}^{k}\left( R\right) $ and $\left( u_{i}^{k}\left(
R\right) \right) ^{\prime }$, both of which are nonnegative, are bounded
above independent of $k$. Using this and the fact that $\left(
u_{i}^{k}\right) ^{\prime }\geq 0$ ($i=1,2$), we note that (\ref{sum})
yields 
\begin{equation*}
\left\{ 
\begin{array}{c}
\left( p-1\right) \left[ \left( u_{1}^{k}\left( r\right) \right) ^{\prime }%
\right] ^{p-1}\left( u_{1}^{k}\right) ^{\prime \prime }\leq a_{1}^{R}\overset%
{2}{\underset{i=1}{\Sigma }}f_{i}\left( \overset{2}{\underset{i=1}{\Sigma }}%
u_{i}^{k},\overset{2}{\underset{i=1}{\Sigma }}u_{i}^{k}\right) \left( 
\overset{2}{\underset{i=1}{\Sigma }}u_{i}^{k}\left( r\right) \right)
^{\prime } \\ 
\left( p-1\right) \left[ \left( u_{2}^{k}\left( r\right) \right) ^{\prime }%
\right] ^{p-1}\left( u_{2}^{k}\right) ^{\prime \prime }\leq a_{2}^{R}\overset%
{2}{\underset{i=1}{\Sigma }}f_{i}\left( \overset{2}{\underset{i=1}{\Sigma }}%
u_{i}^{k},\overset{2}{\underset{i=1}{\Sigma }}u_{i}^{k}\right) \left( 
\overset{2}{\underset{i=1}{\Sigma }}u_{i}^{k}\left( r\right) \right)
^{\prime }%
\end{array}%
\right.
\end{equation*}%
or, equivalently%
\begin{equation*}
\left\{ 
\begin{array}{c}
\frac{p-1}{p}\left\{ \left[ \left( u_{1}^{k}\left( r\right) \right) ^{\prime
}\right] ^{p}\right\} ^{\prime }\leq a_{1}^{R}\overset{2}{\underset{i=1}{%
\Sigma }}f_{i}\left( \overset{2}{\underset{i=1}{\Sigma }}u_{i}^{k},\overset{2%
}{\underset{i=1}{\Sigma }}u_{i}^{k}\right) \left( \overset{2}{\underset{i=1}{%
\Sigma }}u_{i}^{k}\left( r\right) \right) ^{\prime }, \\ 
\frac{p-1}{p}\left\{ \left[ \left( u_{2}^{k}\left( r\right) \right) ^{\prime
}\right] ^{p}\right\} ^{\prime }\leq a_{2}^{R}\overset{2}{\underset{i=1}{%
\Sigma }}f_{i}\left( \overset{2}{\underset{i=1}{\Sigma }}u_{i}^{k},\overset{2%
}{\underset{i=1}{\Sigma }}u_{i}^{k}\right) \left( \overset{2}{\underset{i=1}{%
\Sigma }}u_{i}^{k}\left( r\right) \right) ^{\prime }.%
\end{array}%
\right.
\end{equation*}%
But, this implies%
\begin{equation*}
\left\{ \underset{i=1}{\overset{2}{\Sigma }}\left[ \left( u_{i}^{k}\left(
r\right) \right) ^{\prime }\right] ^{p}\right\} ^{\prime }\leq \frac{p}{p-1}%
\underset{i=1}{\overset{2}{\Sigma }}a_{i}^{R}\overset{2}{\underset{i=1}{%
\Sigma }}f_{i}\left( \overset{2}{\underset{i=1}{\Sigma }}u_{i}^{k}\left(
r\right) ,\overset{2}{\underset{i=1}{\Sigma }}u_{i}^{k}\left( r\right)
\right) \left( \underset{i=1}{\overset{2}{\Sigma }}u_{i}^{k}\left( r\right)
\right) ^{\prime }.
\end{equation*}%
Integrate this equation from $0$ to $r$. We obtain%
\begin{equation}
\underset{i=1}{\overset{2}{\Sigma }}\left[ \left( u_{i}^{k}\left( r\right)
\right) ^{\prime }\right] ^{p}\leq \frac{p}{p-1}\underset{i=1}{\overset{2}{%
\Sigma }}a_{i}^{R}\int_{b}^{\overset{2}{\underset{i=1}{\Sigma }}%
u_{i}^{k}\left( r\right) }\overset{2}{\underset{i=1}{\Sigma }}f_{i}\left(
s,s\right) ds\leq \frac{p}{p-1}\underset{i=1}{\overset{2}{\Sigma }}%
a_{i}^{R}\int_{0}^{\overset{2}{\underset{i=1}{\Sigma }}u_{i}^{k}\left(
r\right) }\overset{2}{\underset{i=1}{\Sigma }}f_{i}\left( s,s\right) ds.
\label{ineq2}
\end{equation}%
Since $p>1$ we know that%
\begin{equation}
\left( a_{1}+a_{2}\right) ^{p}\leq 2^{p-1}\left( a_{1}^{p}+a_{2}^{p}\right)
\label{ineq}
\end{equation}%
for any non-negative constants $a_{i}$ ($i=1,2$). Using this inequality in (%
\ref{ineq2}) we have%
\begin{equation*}
2^{1-p}\left[ \underset{i=1}{\overset{2}{\Sigma }}\left( u_{i}^{k}\left(
r\right) \right) ^{\prime }\right] ^{p}\leq \frac{p}{p-1}\underset{i=1}{%
\overset{2}{\Sigma }}a_{i}^{R}\int_{0}^{\overset{2}{\underset{i=1}{\Sigma }}%
u_{i}^{k}\left( r\right) }\overset{2}{\underset{i=1}{\Sigma }}f_{i}\left(
s,s\right) ds\text{, }0\leq r\leq R,
\end{equation*}%
which yields%
\begin{equation}
\left( \underset{i=1}{\overset{2}{\Sigma }}u_{i}^{k}\left( r\right) \right)
^{\prime }\leq \sqrt[p]{\frac{p2^{p-1}}{p-1}\underset{i=1}{\overset{2}{%
\Sigma }}a_{i}^{R}}\left( \int_{0}^{\overset{2}{\underset{i=1}{\Sigma }}%
u_{i}^{k}\left( r\right) }\overset{2}{\underset{i=1}{\Sigma }}f_{i}\left(
s,s\right) ds\right) ^{1/p}\text{, }0\leq r\leq R.  \label{9}
\end{equation}%
Integrating the above equation between $0$ and $R$, we have 
\begin{equation*}
\int_{b}^{\overset{2}{\underset{i=1}{\Sigma }}u_{i}^{k}\left( R\right) }%
\left[ \int_{0}^{t}\overset{2}{\underset{i=1}{\Sigma }}f_{i}\left(
s,s\right) ds\right] ^{-1/p}dt=I\left( \overset{2}{\underset{i=1}{\Sigma }}%
u_{i}^{k}\left( R\right) \right) -I\left( b\right) \leq \sqrt[p]{\frac{%
p2^{p-1}}{p-1}\left( \underset{i=1}{\overset{2}{\Sigma }}a_{i}^{R}\right) }R%
\text{.}
\end{equation*}%
Since $I$ is a bijection with $I^{-1}$ increasing we obtain%
\begin{equation}
\overset{2}{\underset{i=1}{\Sigma }}u_{i}^{k}\left( R\right) \leq
I^{-1}\left( \sqrt[p]{\frac{p2^{p-1}}{p-1}\left( \underset{i=1}{\overset{2}{%
\Sigma }}a_{i}^{R}\right) }R+I\left( b\right) \right) \text{ for all }r\geq
0.  \label{bound1}
\end{equation}%
By the Keller-Osserman condition (\ref{KO}), we now conclude that $\overset{2%
}{\underset{i=1}{\Sigma }}u_{i}^{k}\left( R\right) $ is uniformly bounded
above independent of $k$ and using this fact in (\ref{9}) shows that the
same is true of $\left( \overset{2}{\underset{i=1}{\Sigma }}u_{i}^{k}\left(
R\right) \right) ^{\prime }$. Thus, the sequences $u_{i}^{k}\left( r\right) $
($i=1,2$) \ are uniformly bounded above independent of $k$ (since $r\leq R$
and $u_{i}^{k}\left( r\right) $ is non-decreasing sequence). Also, we
clearly have $u_{i}^{k}\left( r\right) >0$ for all $r\geq 0$ and so our
sequence is equi-continuous on $\left[ 0,R\right] $ for arbitrary $R>0$. \
Since $u_{i}^{k}\left( r\right) $ ($i=1,2$) is a monotonic, uniformly
bounded, equi-continuous sequence of functions on $\left[ 0,R\right] $ there
exists a function $\left( u_{1},u_{2}\right) \in C\left( \left[ 0,R\right]
\right) \times C\left( \left[ 0,R\right] \right) $ such that $%
u_{i}^{k}\left( r\right) \rightarrow u_{i}^{k}\left( r\right) $ ($i=1,2$)
uniformly. Hence $\left( u_{1},u_{2}\right) $ is a fixed point of (\ref{op})
in $C\left( \left[ 0,R\right] \right) \times C\left( \left[ 0,R\right]
\right) $. Next, we extend this result to show that $S$ has a fixed point in 
$C^{1}\left( \left[ 0,\infty \right) \right) \times C^{1}\left( \left[
0,\infty \right) \right) $. Let $\left\{ u_{i}^{k}\left( r\right) \right\}
_{i=1,2}^{k\geq 1}$ be a sequence of fixed points defined by \ 
\begin{equation}
\left( u_{1}^{k}\left( r\right) ,u_{2}^{k}\left( r\right) \right) =S\left(
u_{1}^{k}\left( r\right) ,u_{2}^{k}\left( r\right) \right) \text{ on }\left[
0,k\right] \text{, }\left( u_{1}^{k}\left( r\right) ,u_{2}^{k}\left(
r\right) \right) \in C\left( \left[ 0,k\right] \right) \times C\left( \left[
0,k\right] \right) ,  \label{c1}
\end{equation}%
for $k=1,2,3,...$ As earlier, we may show that both $u_{1}^{k}\left(
r\right) $ and $u_{2}^{k}\left( r\right) $ are bounded and equi-continuous
on $\left[ 0,1\right] $. Thus by applying the Arzela-Ascoli Theorem to each
sequence separately, we can derive that $\left\{ \left( u_{1}^{k}\left(
r\right) ,u_{2}^{k}\left( r\right) \right) \right\} ^{k\geq 1}$ contains a
convergent subsequence, $\left( u_{1}^{k_{1}^{1}}\left( r\right)
,u_{2}^{k_{2}^{1}}\left( r\right) \right) $, that converges uniformly on $%
\left[ 0,1\right] \times \left[ 0,1\right] $. Let 
\begin{equation*}
\left( u_{1}^{k_{1}^{1}}\left( r\right) ,u_{2}^{k_{2}^{1}}\left( r\right)
\right) \rightarrow \left( u_{1}^{1},u_{2}^{1}\right) \text{ uniformly on }%
\left[ 0,1\right] \times \left[ 0,1\right] \text{ as }k_{1}^{1},k_{2}^{1}%
\rightarrow \infty \text{.}
\end{equation*}%
Likewise, the subsequences $u_{1}^{k_{1}^{1}}\left( r\right) $ and $%
u_{2}^{k_{2}^{1}}\left( r\right) $ are each bounded and equi-continous on $%
\left[ 0,2\right] $ so there exists a subsequence $\left(
u_{1}^{k_{1}^{2}}\left( r\right) ,u_{2}^{k_{2}^{2}}\left( r\right) \right) $
of $\left( u_{1}^{k_{1}^{1}}\left( r\right) ,u_{2}^{k_{2}^{1}}\left(
r\right) \right) $ such that 
\begin{equation*}
\left( u_{1}^{k_{1}^{2}}\left( r\right) ,u_{2}^{k_{2}^{2}}\left( r\right)
\right) \rightarrow \left( u_{1}^{2},u_{2}^{2}\right) \text{ uniformly on }%
\left[ 0,2\right] \times \left[ 0,2\right] \text{ as }k_{1}^{2},k_{2}^{2}%
\rightarrow \infty \text{.}
\end{equation*}%
Notice that 
\begin{equation*}
\left\{ \left( u_{1}^{k_{1}^{2}}\left( r\right) ,u_{2}^{k_{2}^{2}}\left(
r\right) \right) \right\} \subseteq \left\{ \left( u_{1}^{k_{1}^{1}}\left(
r\right) ,u_{2}^{k_{2}^{1}}\left( r\right) \right) \right\} \subseteq
\left\{ \left( u_{1}^{k}\left( r\right) ,u_{2}^{k}\left( r\right) \right)
\right\} _{k\geq 2}^{\infty }
\end{equation*}%
so $\left( u_{1}^{2},u_{2}^{2}\right) =\left( u_{1}^{1},u_{2}^{1}\right) $
on $\left[ 0,1\right] \times \left[ 0,1\right] $. Continuing this line of
reasoning, we obtain a sequence, denoted by $\left\{ \left( u_{1}^{k}\left(
r\right) ,u_{2}^{k}\left( r\right) \right) \right\} $, such that%
\begin{eqnarray*}
\left( u_{1}^{k}\left( r\right) ,u_{2}^{k}\left( r\right) \right) &\in
&C\left( \left[ 0,k\right] \right) \times C\left( \left[ 0,k\right] \right) 
\text{, }k=1,2,.. \\
\left( u_{1}^{k}\left( r\right) ,u_{2}^{k}\left( r\right) \right) &=&\left(
u_{1}^{1}\left( r\right) ,u_{2}^{1}\left( r\right) \right) \text{ for }r\in %
\left[ 0,1\right] \\
\left( u_{1}^{k}\left( r\right) ,u_{2}^{k}\left( r\right) \right) &=&\left(
u_{1}^{2}\left( r\right) ,u_{2}^{2}\left( r\right) \right) \text{ for }r\in %
\left[ 0,2\right] \\
&&... \\
\left( u_{1}^{k}\left( r\right) ,u_{2}^{k}\left( r\right) \right) &=&\left(
u_{1}^{k-1}\left( r\right) ,u_{2}^{k-1}\left( r\right) \right) \text{ for }%
r\in \left[ 0,k-1\right] ,
\end{eqnarray*}%
and these functions are radially symmetric. Therefore $\left(
u_{1}^{k}\left( r\right) ,u_{2}^{k}\left( r\right) \right) $ converges
pointwise to some $\left( u_{1}\left( r\right) ,u_{2}\left( r\right) \right) 
$ which satisfies%
\begin{equation*}
\left( u_{1}\left( r\right) ,u_{2}\left( r\right) \right) =\left(
u_{1}^{k}\left( r\right) ,u_{2}^{k}\left( r\right) \right) \text{ if }0\leq
r\leq k\text{,}
\end{equation*}%
Hence, $\left( u_{1}\left( r\right) ,u_{2}\left( r\right) \right) $ is
radially symmetric. Further, since $\left( u_{1}^{k}\left( r\right)
,u_{2}^{k}\left( r\right) \right) $ is in the form (\ref{c1}) \ we have that 
$\left( u_{1}^{k}\left( r\right) ,u_{2}^{k}\left( r\right) \right) $ is also
equi-continuous.\ \ \ Pointwise convergence and equi-continuity imply
uniform convergence and thus the convergence is uniform on bounded sets.
Thus 
\begin{equation*}
\left( u_{1}\left( r\right) ,u_{2}\left( r\right) \right) \in C^{1}\left( 
\left[ 0,\infty \right) \right) \times C^{1}\left( \left[ 0,\infty \right)
\right)
\end{equation*}%
is a fixed point of (\ref{op}) and a solution to (\ref{11}) with central
value $\left( \frac{b}{2},\frac{b}{2}\right) $. Since $b\geq a>0$ was chosen
arbitrarily, it follows that (\ref{11}) has infinitely many positive entire
solutions and so the first part of our theorem is proved.

\textbf{The proof of i)}\quad Assume that (\ref{5}) holds. Finally, we show
that \textit{any entire positive radial solution }$\left( u_{1},u_{2}\right) 
$\ of system (\ref{11}) is bounded. We choose $R>0$ so that 
\begin{equation*}
r^{\frac{p\left( N-1\right) }{p-1}}\underset{j=1}{\overset{2}{\Sigma }}e^{%
\frac{p}{p-1}\int_{0}^{r}h_{j}\left( t\right) dt}a_{j}\left( r\right) \text{
is non-decreasing for }r\geq R.\ 
\end{equation*}%
Multiplying%
\begin{equation}
\left\{ 
\begin{array}{c}
\left( p-1\right) \left[ \left( u_{1}\left( r\right) \right) ^{\prime }%
\right] ^{p-1}\left( u_{1}\right) ^{\prime \prime }+\left( \frac{N-1}{r}%
+h_{1}\left( r\right) \right) \left[ \left( u_{1}\left( r\right) \right)
^{\prime }\right] ^{p}\leq a_{1}\left( r\right) \overset{2}{\underset{i=1}{%
\Sigma }}f_{i}\left( \overset{2}{\underset{i=1}{\Sigma }}u_{i},\overset{2}{%
\underset{i=1}{\Sigma }}u_{i}\right) \left( \overset{2}{\underset{i=1}{%
\Sigma }}u_{i}\left( r\right) \right) ^{\prime }, \\ 
\left( p-1\right) \left[ \left( u_{2}\left( r\right) \right) ^{\prime }%
\right] ^{p-1}\left( u_{2}\right) ^{\prime \prime }+\left( \frac{N-1}{r}%
+h_{2}\left( r\right) \right) \left[ \left( u_{2}\left( r\right) \right)
^{\prime }\right] ^{p}\leq a_{2}\left( r\right) \overset{2}{\underset{i=1}{%
\Sigma }}f_{i}\left( \overset{2}{\underset{i=1}{\Sigma }}u_{i},\overset{2}{%
\underset{i=1}{\Sigma }}u_{i}\right) \left( \overset{2}{\underset{i=1}{%
\Sigma }}u_{i}\left( r\right) \right) ^{\prime }.%
\end{array}%
\right.  \notag
\end{equation}%
each line of this system by 
\begin{equation*}
\frac{p}{p-1}r^{\frac{p\left( N-1\right) }{p-1}}e^{\frac{p}{p-1}%
\int_{0}^{r}h_{i}(t)dt}\text{, (}i=1,2\text{)}
\end{equation*}%
($i$\ represent the equation of the system that will be multiplied) and
summing we have%
\begin{equation*}
\left\{ r^{\frac{p\left( N-1\right) }{p-1}}\overset{2}{\underset{i=1}{\Sigma 
}}e^{\frac{p}{p-1}\int_{0}^{r}h_{i}(t)dt}\left[ \left( u_{i}\right) ^{\prime
}\right] ^{p}\right\} ^{\prime }\leq \text{ }\frac{pr^{\frac{p\left(
N-1\right) }{p-1}}}{p-1}\overset{2}{\underset{i=1}{\Sigma }}e^{\frac{p}{p-1}%
\int_{0}^{r}h_{i}(t)dt}a_{i}\left( r\right) \overset{2}{\underset{i=1}{%
\Sigma }}f_{i}\left( \overset{2}{\underset{i=1}{\Sigma }}u_{i},\overset{2}{%
\underset{i=1}{\Sigma }}u_{i}\right) \left( \overset{2}{\underset{i=1}{%
\Sigma }}u_{i}\right) ^{\prime }.
\end{equation*}%
and integrating gives%
\begin{eqnarray}
&&\int_{R}^{r}\left\{ s^{\frac{p\left( N-1\right) }{p-1}}\overset{2}{%
\underset{i=1}{\Sigma }}\left[ e^{\frac{1}{p-1}\int_{0}^{r}h_{i}(t)dt}\left(
u_{i}\left( s\right) \right) ^{\prime }\right] ^{p}\right\} ^{\prime }ds 
\notag \\
&\leq &\int_{R}^{r}\frac{p}{p-1}s^{\frac{p\left( N-1\right) }{p-1}}\overset{2%
}{\underset{i=1}{\Sigma }}e^{\frac{p}{p-1}\int_{0}^{s}h_{i}(t)dt}a_{i}\left(
s\right) \overset{2}{\underset{i=1}{\Sigma }}f_{i}\left( \overset{2}{%
\underset{i=1}{\Sigma }}u_{i}\left( s\right) ,\overset{2}{\underset{i=1}{%
\Sigma }}u_{i}\left( s\right) \right) \left( \overset{2}{\underset{i=1}{%
\Sigma }}u_{i}\right) ^{\prime }ds.  \label{s2}
\end{eqnarray}%
Hence, using (\ref{ineq}) in (\ref{s2}) it gives 
\begin{eqnarray*}
&&r^{\frac{p\left( N-1\right) }{p-1}}2^{1-p}\left[ \left( \overset{2}{%
\underset{i=1}{\Sigma }}e^{\frac{1}{p-1}\int_{0}^{r}h_{i}(t)dt}u_{i}\left(
r\right) \right) ^{\prime }\right] ^{p}-R^{\frac{p\left( N-1\right) }{p-1}}%
\overset{2}{\underset{i=1}{\Sigma }}\left[ e^{\frac{1}{p-1}%
\int_{0}^{r}h_{i}(t)dt}\left( u_{i}\left( R\right) \right) ^{\prime }\right]
^{p} \\
&\leq &\int_{R}^{r}\frac{p}{p-1}s^{\frac{p\left( N-1\right) }{p-1}}\overset{2%
}{\underset{i=1}{\Sigma }}e^{\frac{p}{p-1}\int_{0}^{s}h_{i}(t)dt}a_{i}\left(
s\right) \overset{2}{\underset{i=1}{\Sigma }}f_{i}\left( \overset{2}{%
\underset{i=1}{\Sigma }}u_{i}\left( s\right) ,\overset{2}{\underset{i=1}{%
\Sigma }}u_{i}\left( s\right) \right) \left( \overset{2}{\underset{i=1}{%
\Sigma }}u_{i}\right) ^{\prime }ds
\end{eqnarray*}%
and thus%
\begin{eqnarray*}
&&r^{\frac{p\left( N-1\right) }{p-1}}\left[ \left( \overset{2}{\underset{i=1}%
{\Sigma }}e^{\frac{1}{p-1}\int_{0}^{r}h_{i}(t)dt}u_{i}\left( r\right)
\right) ^{\prime }\right] ^{p}\leq R^{\frac{p\left( N-1\right) }{p-1}}2^{p-1}%
\overset{2}{\underset{i=1}{\Sigma }}\left[ e^{\frac{1}{p-1}%
\int_{0}^{r}h_{i}(t)dt}\left( u_{i}\left( R\right) \right) ^{\prime }\right]
^{p}+ \\
&&+\int_{R}^{r}\frac{p2^{p-1}s^{\frac{p\left( N-1\right) }{p-1}}}{p-1}%
\overset{2}{\underset{i=1}{\Sigma }}e^{\frac{p}{p-1}%
\int_{0}^{s}h_{i}(t)dt}a_{i}\left( s\right) \overset{2}{\underset{i=1}{%
\Sigma }}f_{i}\left( \overset{2}{\underset{i=1}{\Sigma }}u_{i}\left(
s\right) ,\overset{2}{\underset{i=1}{\Sigma }}u_{i}\left( s\right) \right)
\left( \overset{2}{\underset{i=1}{\Sigma }}u_{i}\right) ^{\prime }ds.
\end{eqnarray*}%
for $r\geq R$. Noting that, by the monotonicity of $s^{\frac{p\left(
N-1\right) }{p-1}}\overset{2}{\underset{i=1}{\Sigma }}e^{\frac{p}{p-1}%
\int_{0}^{s}h_{i}(t)dt}a_{i}\left( s\right) $ for $r\geq s\geq R$, we get 
\begin{equation*}
r^{^{\frac{p\left( N-1\right) }{p-1}}}\left[ \left( \overset{2}{\underset{i=1%
}{\Sigma }}e^{\frac{1}{p-1}\int_{0}^{r}h_{i}(t)dt}u_{i}\left( r\right)
\right) ^{\prime }\right] ^{p}\leq C+\frac{p2^{p-1}}{p-1}r^{^{\frac{p\left(
N-1\right) }{p-1}}}\overset{2}{\underset{i=1}{\Sigma }}e^{\frac{p}{p-1}%
\int_{0}^{r}h_{i}(t)dt}a_{i}\left( r\right) F\left( \overset{2}{\underset{i=1%
}{\Sigma }}u_{i}\left( r\right) \right) ,
\end{equation*}%
where 
\begin{equation*}
C=R^{^{\frac{p\left( N-1\right) }{p-1}}}2^{p-1}\overset{2}{\underset{i=1}{%
\Sigma }}\left[ e^{\frac{1}{p-1}\int_{0}^{R}h_{i}(t)dt}\left( u_{i}\left(
R\right) \right) ^{\prime }\right] ^{p},
\end{equation*}%
which yields%
\begin{equation}
\left( \overset{2}{\underset{i=1}{\Sigma }}e^{\frac{1}{p-1}%
\int_{0}^{r}h_{i}(t)dt}u_{i}\right) ^{\prime }\leq \left[ Cr^{^{\frac{%
p\left( 1-N\right) }{p-1}}}+\frac{p2^{p-1}}{p-1}\overset{2}{\underset{i=1}{%
\Sigma }}e^{\frac{p}{p-1}\int_{0}^{r}h_{i}(t)dt}a_{i}\left( r\right) F\left( 
\overset{2}{\underset{i=1}{\Sigma }}u_{i}\left( r\right) \right) \right]
^{1/p}.  \label{100}
\end{equation}%
Since $\left( 1/p\right) <1$ we know that%
\begin{equation*}
\left( b_{1}+b_{2}\right) ^{1/p}\leq b_{1}^{1/p}+b_{2}^{1/p}
\end{equation*}%
for any non-negative constants $b_{i}$ ($i=1,2$). Therefore, by applying
this inequality in (\ref{100}) we get%
\begin{eqnarray*}
\left( \overset{2}{\underset{i=1}{\Sigma }}e^{\frac{1}{p-1}%
\int_{0}^{r}h_{i}(t)dt}u_{i}\right) ^{\prime } &\leq &\sqrt[p]{C}r^{\left(
1-N\right) /(p-1)}+\sqrt[p]{\frac{p2^{p-1}}{p-1}\overset{2}{\underset{i=1}{%
\Sigma }}e^{\frac{p}{p-1}\int_{0}^{r}h_{i}(t)dt}a_{i}\left( r\right) }\left[
F\left( \overset{2}{\underset{i=1}{\Sigma }}u_{i}\left( r\right) \right) %
\right] ^{1/p} \\
&\leq &\sqrt[p]{C}r^{\left( 1-N\right) /(p-1)}+\sqrt[p]{\frac{p2^{p-1}}{p-1}%
\overset{2}{\underset{i=1}{\Sigma }}e^{\frac{p}{p-1}%
\int_{0}^{r}h_{i}(t)dt}a_{i}\left( r\right) }\left[ F\left( \overset{2}{%
\underset{i=1}{\Sigma }}e^{\frac{1}{p-1}\int_{0}^{r}h_{i}(t)dt}u_{i}\left(
r\right) \right) \right] ^{1/p}.
\end{eqnarray*}%
Integrating the above inequality, we get%
\begin{eqnarray}
&&\frac{d}{dr}\int_{\overset{2}{\underset{i=1}{\Sigma }}e^{\frac{1}{p-1}%
\int_{0}^{r}h_{i}(t)dt}u_{i}\left( R\right) }^{\overset{2}{\underset{i=1}{%
\Sigma }}e^{\frac{1}{p-1}\int_{0}^{r}h_{i}(t)dt}u_{i}\left( r\right) }\left[
F\left( t\right) \right] ^{-1/p}dt  \label{sis3} \\
&\leq &\sqrt[p]{C}r^{\left( 1-N\right) /\left( p-1\right) }\left[ F\left( 
\overset{2}{\underset{i=1}{\Sigma }}e^{\frac{1}{p-1}%
\int_{0}^{r}h_{i}(t)dt}u_{i}\left( r\right) \right) \right] ^{-1/p}+\left( 
\frac{p2^{p-1}}{p-1}\overset{2}{\underset{i=1}{\Sigma }}e^{\frac{p}{p-1}%
\int_{0}^{r}h_{i}(t)dt}a_{i}\left( r\right) \right) ^{1/p}.  \notag
\end{eqnarray}%
Integrating (\ref{sis3}) and using the fact that%
\begin{eqnarray*}
\left( \overset{2}{\underset{i=1}{\Sigma }}e^{\frac{p}{p-1}%
\int_{0}^{r}h_{i}(t)dt}\varphi _{i}\left( s\right) \right) ^{1/p} &=&\left(
s^{p\left( 1+\varepsilon \right) /2}\overset{2}{\underset{i=1}{\Sigma }}e^{%
\frac{p}{p-1}\int_{0}^{r}h_{i}(t)dt}a_{i}\left( s\right) s^{-p\left(
1+\varepsilon \right) /2}\right) ^{1/p} \\
&\leq &\left( \frac{1}{2}\right) ^{1/p}\left[ s^{1+\varepsilon }\left( 
\overset{2}{\underset{i=1}{\Sigma }}e^{\frac{p}{p-1}%
\int_{0}^{r}h_{i}(t)dt}a_{i}\left( r\right) \right) ^{2/p}+s^{-1-\varepsilon
}\right] ,
\end{eqnarray*}%
for each $\varepsilon >0$, we have%
\begin{eqnarray}
&&\int_{\overset{2}{\underset{i=1}{\Sigma }}e^{\frac{1}{p-1}%
\int_{0}^{r}h_{i}(t)dt}u_{i}\left( R\right) }^{\overset{2}{\underset{i=1}{%
\Sigma }}e^{\frac{1}{p-1}\int_{0}^{r}h_{i}(t)dt}u_{i}\left( r\right) }\left[
F\left( t\right) \right] ^{-1/p}dt  \notag \\
&\leq &\sqrt[p]{C}\int_{R}^{r}t^{\frac{1-N}{p-1}}\left[ F\left( \overset{2}{%
\underset{i=1}{\Sigma }}e^{\frac{1}{p-1}\int_{0}^{t}h_{i}(t)dt}u_{i}\left(
t\right) \right) \right] ^{-1/p}dt  \notag \\
&&+\left( \frac{1}{2}\right) ^{1/p}\sqrt[p]{\frac{p2^{p-1}}{p-1}}\left[
\int_{R}^{r}t^{1+\varepsilon }\left( \overset{2}{\underset{i=1}{\Sigma }}e^{%
\frac{p}{p-1}\int_{0}^{t}h_{i}(t)dt}a_{i}\left( t\right) \right)
^{2/p}dt+\int_{R}^{r}t^{-1-\varepsilon }dt\right]  \notag \\
&\leq &\sqrt[p]{C}\left[ F\left( \overset{2}{\underset{i=1}{\Sigma }}e^{%
\frac{1}{p-1}\int_{0}^{R}h_{i}(t)dt}u_{i}\left( R\right) \right) \right]
^{-1/p}\frac{p-1}{p-N}R^{\frac{p-N}{p-1}}  \notag \\
&&+\left( \frac{1}{2}\right) ^{1/p}\sqrt[p]{\frac{p2^{p-1}}{p-1}}\left[
\int_{R}^{r}t^{1+\varepsilon }\left( \overset{2}{\underset{i=1}{\Sigma }}e^{%
\frac{p}{p-1}\int_{0}^{t}h_{i}(t)dt}a_{i}\left( t\right) \right) ^{2/p}dt+%
\frac{1}{\varepsilon R^{\varepsilon }}\right] \text{.}  \label{111}
\end{eqnarray}%
Since the right side of this inequality is bounded (note that $u_{i}\left(
t\right) \geq b/2$), so is the left side and hence, in light of Keller
Osserman condition, the sequence $\overset{2}{\underset{i=1}{\Sigma }}e^{%
\frac{1}{p-1}\int_{0}^{r}h_{i}(t)dt}u_{i}\left( r\right) $ is bounded and so 
$\overset{2}{\underset{i=1}{\Sigma }}u_{i}\left( r\right) $ that implies
finally $u_{i}\left( r\right) $ ($i=1,2$) is a bounded function. \ Thus, for
every $x\in \mathbb{R}^{N}$ $\left( u_{1}\left( \left\vert x\right\vert
\right) ,u_{2}\left( \left\vert x\right\vert \right) \right) $ is a positive
bounded solution of (\ref{11}).

\textbf{The proof of ii)\quad }Suppose that $a_{i}$ ($i=1,2$) satisfies (\ref%
{12}). Now, let $\left( u_{1},u_{2}\right) $ be any positive entire radial
solution of (\ref{11}) determined in the first step of the proof. Since $%
u_{i}$ ($i=1,2$) is positive for all $R>0$ we have $u_{i}\left( R\right) >0$%
. Since $u_{i}^{\prime }\geq 0$, we get $u_{i}\left( r\right) \geq
u_{i}\left( R\right) $ for $r\geq R$ and thus from%
\begin{equation*}
u_{i}\left( r\right) =u_{i}\left( 0\right) +\int_{0}^{r}\frac{%
e^{-\int_{0}^{t}h_{i}\left( s\right) ds}}{t^{N-1}}\left(
\int_{0}^{t}s^{N-1}e^{\int_{0}^{s}h_{i}\left( s\right) ds}a_{i}\left(
s\right) f_{i}\left( u_{1}\left( s\right) ,u_{2}\left( s\right) \right)
ds\right) ^{1/\left( p-1\right) }dt,
\end{equation*}%
we obtain

\begin{equation*}
\left\{ 
\begin{array}{l}
u_{i}\left( r\right) =u_{i}\left( 0\right) +\int_{0}^{r}\left( \frac{%
e^{-\int_{0}^{t}h_{i}\left( s\right) ds}}{t^{N-1}}\int_{0}^{t}s^{N-1}e^{%
\int_{0}^{t}h_{i}\left( s\right) ds}a_{i}\left( s\right) f_{i}\left(
u_{1}\left( s\right) ,u_{2}\left( s\right) \right) ds\right) ^{1/\left(
p-1\right) }dt\text{ }\geq u_{i}\left( R\right) \\ 
\text{ \ \ \ \ \ \ }+f_{i}^{1/\left( p-1\right) }\left( u_{1}\left( R\right)
,u_{2}\left( R\right) \right) \int_{R}^{r}\left( \frac{e^{-\int_{0}^{t}h_{i}%
\left( s\right) ds}}{t^{N-1}}\int_{R}^{t}s^{N-1}e^{\int_{0}^{s}h_{i}\left(
s\right) ds}a_{i}\left( s\right) ds\right) ^{1/\left( p-1\right)
}dt\rightarrow \infty \text{ as }r\rightarrow \infty , \\ 
\text{for all }i=1,2.%
\end{array}%
\right.
\end{equation*}%
and the proof is complete.

From the above proof and the work \cite{CD2} we can easy obtain the
following remarks.

\begin{remark}
Make the same assumptions as in Theorem \ref{1} on $a_{j}$, $h_{j}$ and $%
f_{j}$ excepting i)-ii). If, on the other hand, $a_{j}$\ satisfies%
\begin{equation}
\left( \frac{1}{N}\right) ^{1/\left( p-1\right) }\int_{0}^{\infty }\left(
e^{-\int_{0}^{t}h_{j}\left( s\right) ds}ta_{j}\left( t\right) \right)
^{1/\left( p-1\right) }dt=\infty ,\text{ }j=1,2  \label{5b}
\end{equation}%
then system (\ref{non}) has no nonnegative nontrivial entire bounded radial
solution on $\mathbb{R}^{N}$\textit{.}\bigskip
\end{remark}

\begin{remark}
(see and \cite{CD2} for the proof) Make the same assumptions as in Theorem %
\ref{1} on $a_{j}$, $h_{j}$ and $f_{j}$ excepting i)-ii).\textit{\ If (\ref%
{non}) has a nonnegative entire large solution, then }$a_{j}$\textit{\ }($%
j=1,2$) \textit{satisfy}%
\begin{equation}
\int_{0}^{\infty }r^{1+\varepsilon }\left( \underset{j=1}{\overset{2}{\Sigma 
}}e^{\frac{p}{p-1}\int_{0}^{t}h_{j}(t)dt}a_{j}\left( t\right) \right)
^{2/p}dr=\infty ,  \label{13}
\end{equation}%
\textit{for every }$\varepsilon >0$\textit{.}
\end{remark}

\begin{remark}
As we have observed in the article \cite{CD2} the above proofs can be
adopted to obtain the same results for systems with indefinite number of
equations.
\end{remark}

\bigskip

\bigskip

\end{document}